\documentclass[11pt]{article}
\usepackage{amsmath}
\usepackage{amssymb}
\usepackage{amscd}
\usepackage[all]{xy}
\usepackage[T1]{fontenc}
\usepackage{epsfig}
\def\B{\mathcal B}
\def\qsl{\mathcal O(SL_q(2))}
\def\N{\mathbb N}
\def\C{\mathbb C}

\newtheorem{theo}{Theorem}[section]
\newtheorem{prop}[theo]{Proposition}
\newtheorem{coro}[theo]{Corollary}
\newtheorem{lemm}[theo]{Lemma}
\newtheorem{defi}[theo]{Definition}
\newtheorem{rem}[theo]{Remark}
\newtheorem{ex}[theo]{Example}

\title{\textbf{\textsc{Hopf-Galois Systems}}} 
\author{\textsc{Julien Bichon}}
\date{{\small \textsl{Laboratoire de Math\'ematiques Appliqu\'ees,
Universit\'e de Pau et des Pays de l'Adour, \\
IPRA, Avenue de l'universit\'e,
64000 Pau, France.}
E-mail: Julien.Bichon@univ-pau.fr}}

\makeatletter
\renewcommand{\@makefnmark}{}
\makeatother

\begin{document}

\maketitle

\begin{abstract}
We introduce the concept of Hopf-Galois system,
a reformulation of the notion of Galois extension of the base field
for a Hopf algebra. The main feature of our definition is a generalization
of the antipode of an ordinary Hopf algebra. 
We present several examples which indicate that although our axiomatic is
slightly more complicated than the classical one, it is also
more natural and easier to handle with.
The main application of Hopf-Galois systems is the construction of
monoidal equivalences between comodule categories.
\end{abstract}


\section*{Introduction}

We introduce the concept of Hopf-Galois system, a reformulation of the
notion of Galois extension of the base field for a Hopf algebra.
Our motivation for such a definition is to provide a natural 
way to construct
monoidal equivalences between categories of comodules over Hopf algebras.

\smallskip

Let $A$ and $B$ be Hopf algebras
(over a field $k$). Recall \cite{[Sc1]} that a non-zero algebra
$Z$ is said to be an $A$-$B$-biGalois extension if 
$Z$ is an $A$-$B$-bicomodule algebra such that
two linear maps $\kappa_l : Z \otimes Z \longrightarrow A \otimes Z$
and $\kappa_r : Z \otimes Z \longrightarrow Z  \otimes B$
are bijective (see Section 1). A useful theorem of
P. Schauenburg \cite{[Sc1]} brings interest for biGalois extensions:
the comodule categories over $A$ and $B$ are monoidally equivalent
if and only if there exists an $A$-$B$-biGalois extension.

In \cite{[Bi1],[Bi3]} we constructed some examples of biGalois extensions.
It is more or less clear in these papers that, when checking that
the maps $\kappa_l$ and $\kappa_r$ are bijective, one uses the
concept of Hopf-Galois system introduced in the present paper.
In fact the construction of a Hopf-Galois system seems to be the easiest 
and most natural way to get a biGalois extension. Here is a related
and well-known situation.
Consider a bialgebra $A$. Then
$A$ is a Hopf algebra if only if the the map
\begin{equation*}
\begin{CD}
\kappa_l : A \otimes A @>\Delta \otimes 1_A>>
A \otimes A \otimes A @>1_A \otimes m>> A \otimes A
\end{CD}
\end{equation*}
is bijective (this is well-know, e.g. to multiplier Hopf algebraists
\cite{[VD]}). In concrete examples, it is much more desirable
to require the existence of an antipode, although the axiomatic
is slightly more involved. We adopt the same philosophy for Galois
extensions.

\medskip

A Hopf-Galois system consists of four non-zero algebras
$(A,B,Z,T)$. The algebras $A$ and $B$ are bialgebras, and $Z$ is assumed to be 
an $A$-$B$-bicomodule algebra. There are also algebra morphisms
$\gamma : A \longrightarrow Z \otimes T$ and 
$\delta : B \longrightarrow T \otimes Z$ with some associativity conditions.
Finally there is a linear map $S : T  \longrightarrow Z$
playing the role of an antipode. 
In fact the first axioms are closely related to the ones of a set
of pre-equivalence data of M. Takeuchi (Definition 2.3 in 
\cite{[Ta2]}), the main new feature being the generalized antipode.
See Section 1 for the details.
The easiest way to understand the axioms is to see a Hopf-Galois system
as the dual object of a groupoid with two objects, with some
structures forgotten. In fact a Hopf-Galois system is always associated
to a pair of objects in the groupoid of fibre functors over the comodules
of a Hopf algebra.    
We show (Theorem 1.2) that if  $(A,B,Z,T)$ is a Hopf-Galois system, then
$Z$ is an $A$-$B$-biGalois extension. Conversely, starting from a Galois 
extension, it is possible to reconstruct a Hopf-Galois system.

The axiomatic of a Hopf-Galois system is more complicated than the one 
of an $A$-Galois extension or of an $A$-$B$-biGalois extension. 
But in a completely parallel way to the observation 
concerning bialgebras and Hopf algebras, it is also
very natural and easy to handle with when dealing 
with concrete examples. In fact when
one suspects an algebra to be an $A$-Galois extension, it is not 
difficult to guess what the whole Hopf-Galois system will be. 
We present several examples, that, we hope, will convince the reader.

\smallskip

It is quite possible that the axiomatic of Hopf-Galois systems
was already known to some experts.
It is more or less implicit in \cite{[Sc1]}, with the notation 
$Z^{-1}$ for the fourth algebra.
On the other hand we feel that it should be useful to have 
written down the axiomatic completely, especially in view of applications
in representation theory.

\medskip

Our work is organized as follows.
In Section 1 we give the full definition of a Hopf-Galois system, and 
we show that such a system always gives rise to a biGalois extension.
We also discuss the reconstruction of a Hopf-Galois system 
from a Galois extension.
Since we have decided in this paper to concentrate on examples and 
applications rather than on theoretical aspects of Hopf-Galois systems,
the other sections are devoted to examples. We study the following families of 
examples.

\noindent
$\bullet$ Hopf-Galois systems associated to $2$-cocycles,
with emphasis on the function algebra on the symmetric group, and some
examples generalizing those of  
\cite{[Bi1]} are presented.

\noindent
$\bullet$ Hopf-Galois systems for the Hopf algebra 
of a non-degenerate bilinear form \cite{[DVL], [Bi3]}.

\noindent
$\bullet$ Hopf-Galois systems for universal cosovereign Hopf algebras
\cite{[Bi2]}. This gives some improvements on the known results \cite{[Ba]}
concerning the corepresentation theory of these Hopf algebras.

\noindent
$\bullet$ Hopf-Galois systems for free Hopf algebras generated
by matrix coalgebras \cite{[T]}.

\medskip

\noindent
\textbf{Notations and conventions}.
Throughout this paper $k$ denotes a commutative field.
The reader is assumed to be familiar with Hopf algebras, their modules,
comodules, comodule algebras \cite{[Mo]}. We also
assume familiarity with monoidal categories, monoidal functors
\cite{[Ka], [JS]}. The monoidal category of comodules 
(resp. finite-dimensional comodules) over a 
bialgebra $A$ is denoted by Comod($A$) (resp. Comod$_{\rm f}$($A$)).
The monoidal category of finite-dimensional vector spaces over $k$
is denoted by ${\rm Vect}_{\rm f}(k)$.

\section{Definition and basic results}

We present the formal definition of a Hopf-Galois system.
We work in the monoidal category of vector spaces over $k$, but it is clear 
that our definition still makes sense in any braided monoidal category,
and that Theorem 1.2 is valid in such a category. 
Let us first recall the language of Galois extensions for
Hopf algebras (see \cite{[Mo]} for a general perspective). 

Let $A$ be a Hopf algebra. A left $A$-Galois extension (of $k$)
is a non-zero left $A$-comodule algebra $Z$ such that the linear map
$\kappa_l$ defined by the composition
\begin{equation*}
\begin{CD}
\kappa_l : Z \otimes Z @>\alpha \otimes 1_Z>>
A \otimes Z \otimes Z @>1_A \otimes m_Z>> A \otimes Z
\end{CD}
\end{equation*}
where $\alpha$ is the coaction of $A$ and $m_Z$ is the multiplication
of $Z$, is bijective.

Similarly, a right $A$-Galois extension 
is a non-zero right $A$-comodule algebra $Z$ such that the linear map
$\kappa_r$ defined by the composition
\begin{equation*}
\begin{CD}
\kappa_r : Z \otimes Z @>1_Z \otimes \beta>>
Z \otimes Z \otimes A @>m_Z \otimes 1_A>> Z \otimes A
\end{CD}
\end{equation*}
where $\beta$ is the coaction of $A$, is bijective.

Let $A$ and $B$ be Hopf algebras. An algebra $Z$ 
is said to be an $A$-$B$-bigalois extension \cite{[Sc1]} if $Z$
is both a left $A$-Galois extension and a right $B$-Galois extension,
and if $Z$ is an $A$-$B$-bicomodule.

\begin{defi} A Hopf-Galois system consists of four non-zero algebras
$(A,B,Z,T)$, with the following axioms.

\noindent
${\mathsf{(HG1)}}$ The algebras $A$ and $B$ are bialgebras.

\smallskip

\noindent
${\mathsf{(HG2)}}$ The algebra $Z$ is an $A$-$B$-bicomodule algebra.

\smallskip

\noindent
${\mathsf{(HG3)}}$ There are algebra morphisms
$\gamma : A \longrightarrow  Z  \otimes T$ and
$\delta : B \longrightarrow T \otimes Z$ such that 
the following diagrams commute:
\begin{equation*}
\begin{CD}
Z 
@>\alpha>> A \otimes Z \\
@V\beta VV
@V\gamma \otimes 1_ZVV\\
Z \otimes B  @>1_Z \otimes \delta >> 
Z\otimes T \otimes Z 
\end{CD} \quad
\begin{CD}
A 
@>\Delta_A>> A \otimes A \\
@V\gamma VV
@V1_A \otimes \gamma VV\\
Z \otimes T  @>\alpha \otimes 1_T >> 
A\otimes Z \otimes T 
\end{CD} \quad
\begin{CD}
B 
@>\Delta_B>> B \otimes B\\
@V\delta VV
@V\delta \otimes 1_BVV\\
T\otimes Z @>1_T \otimes \beta>> 
T\otimes Z \otimes B 
\end{CD}
\end{equation*}

\noindent
${\mathsf{(HG4)}}$ There is a linear map 
$S : T \longrightarrow Z$ such that the following diagrams commute:
$$\xymatrix{
A \ar[r]^{\varepsilon_A} \ar[d]^{\gamma} &
k \ar[r]^{u_Z}  & Z \\
Z \otimes T \ar[rr]^{1_Z \otimes S} && Z \otimes Z \ar[u]^{m_Z}} \quad \quad
\xymatrix{
B \ar[r]^{\varepsilon_B} \ar[d]^{\delta} &
k \ar[r]^{u_Z}  & Z \\
T \otimes Z \ar[rr]^{S \otimes 1_Z} && Z \otimes Z \ar[u]^{m_Z}}
$$
\end{defi}

When $A= B = Z= T$ and $\alpha = \beta = \gamma = \delta$, we just 
have the axioms of a Hopf algebra, the linear map $S$ 
being the antipode. The axiom \textsf{HG3} is 
very close from the axioms of a set of pre-equivalence data 
of M. Takeuchi (\cite{[Ta2]}), but we do not require a $B$-$A$-bicomodule
structure on $T$.
We have already mentioned the fact that the easiest way to understand 
the axioms of a Hopf-Galois system is to see it as the dual object
of a groupoid with two objects, where some structures would have
been forgotten. In fact, we have only included the axioms needed
to prove the following result. 

\begin{theo}
Let $(A,B,Z,T)$ be a Hopf-Galois system. Then 
$Z$ is an $A$-$B$-biGalois extension.
\end{theo}

\noindent
\textbf{Proof}. Let $\eta_l : A \otimes Z \longrightarrow Z \otimes Z$
be the morphism defined by 
$$\eta_l = (1_Z \otimes m_Z) \circ (1_Z \otimes S \otimes 1_Z)
\circ (\gamma \otimes 1_Z).$$ 
We show that $\eta_l$ is an inverse for $\kappa_l$.
we have
\begin{align*}
\eta_l \circ \kappa_l & =
(1_Z \otimes m_Z) \circ (1_Z \otimes S \otimes 1_Z)
\circ (\gamma \otimes 1_Z) \circ (1_A \otimes m_Z) \circ
(\alpha \otimes 1_Z) =\\
(1_Z \otimes & m_Z) \circ (1_Z \otimes 1_Z \otimes m_Z) \circ
(1_Z \otimes S \otimes 1_Z \otimes 1_Z) \circ (\gamma \otimes 1_Z \otimes 1_Z)
\circ  (\alpha \otimes 1_Z) = \\
 (1_Z \otimes & m_Z) \circ (1_Z \otimes m_Z \otimes 1_Z) \circ
(1_Z \otimes S \otimes 1_Z \otimes 1_Z) \circ (1_Z \otimes \delta \otimes 1_Z)
\circ (\beta \otimes 1_Z) = \\
(1_Z \otimes & m_Z) \circ (1_Z \otimes (u_Z \circ \varepsilon_B) \otimes
1_Z) \circ (\beta \otimes 1_Z) = 1_{Z \otimes Z}.
\end{align*}
We also have
\begin{align*}
\kappa_l \circ \eta_l & = 
(1_A \otimes m_Z) \circ (\alpha \otimes 1_Z) \circ
(1_Z \otimes m_Z) \circ (1_Z \otimes S \otimes 1_Z)
\circ (\gamma \otimes 1_Z) = \\
(1_A \otimes & m_Z) \circ (1_A \otimes 1_Z \otimes m_Z) \circ
(1_A \otimes 1_Z \otimes S \otimes 1_Z) \circ 
(\alpha \otimes 1_T \otimes 1_Z) \circ (\gamma \otimes 1_Z) = \\
(1_A \otimes & m_Z) \circ (1_A \otimes m_Z \otimes 1_Z) \circ
(1_A \otimes 1_Z \otimes S \otimes 1_Z) \circ 
(1_A \otimes \gamma \otimes 1_Z) \circ (\Delta_A \otimes 1_Z) = \\
(1_A \otimes & m_Z) \circ (1_A \otimes (u_Z \circ \varepsilon_A) \otimes 1_Z)
\circ (\Delta_A \otimes 1_Z) = 1_{A \otimes Z}.
\end{align*}
This proves that $\kappa_l$ is an isomorphism. Similarly,
we define a morphism $\eta_r : Z \otimes B \longrightarrow Z \otimes Z$
by 
$$\eta_r = (m_Z \otimes 1_Z) \circ (1_Z \otimes S \otimes 1_Z)
\circ (1_Z \otimes \delta),$$
and show that $\eta_r$ is an inverse for $\kappa_r$. We have
\begin{align*}
\kappa_r \circ \eta_r & =
(m_Z \otimes 1_B) \circ (1_Z \otimes \beta) \circ
(m_Z \otimes 1_Z) \circ (1_Z \otimes S \otimes 1_Z)
\circ (1_Z \otimes \delta) = \\
(m_Z \otimes & 1_B) \circ (m_Z \otimes 1_Z \otimes 1_B) \circ
(1_Z \otimes S \otimes 1_Z \otimes 1_B) \circ 
(1_Z \otimes 1_T \otimes \beta) \circ (1_Z \otimes \delta) = \\
(m_Z \otimes & 1_B) \circ (1_Z \otimes m_Z \otimes 1_B) \circ
(1_Z \otimes S \otimes 1_Z \otimes 1_B) \circ 
(1_Z \otimes \delta \otimes 1_B) \circ (1_Z \otimes \Delta_B) = \\
(m_Z \otimes & 1_B) \circ 
(1_Z \otimes (u_B \circ \varepsilon_B) \otimes 1_B) 
\circ (1_Z \otimes \Delta_B) = 1_{Z \otimes B}.
\end{align*}
We also have
\begin{align*}
\eta_r \circ \kappa_r & =
(m_Z \otimes 1_Z) \circ (1_Z \otimes S \otimes 1_Z) \circ (1_Z \otimes \delta)
\circ (m_Z \otimes 1_B) \circ (1_Z \otimes \beta) = \\
(m_Z \otimes & 1_Z) \circ (m_Z \otimes 1_Z \otimes 1_Z)
\circ (1_Z \otimes 1_Z \otimes S \otimes 1_Z) \circ
(1_Z \otimes 1_Z \otimes \delta) \circ (1_Z \otimes \beta) = \\
(m_Z \otimes & 1_Z) \circ (1_Z \otimes m_Z \otimes 1_Z)
\circ (1_Z \otimes 1_Z \otimes S \otimes 1_Z) \circ
(1_Z \otimes \gamma \otimes 1_Z) \circ (1_Z \otimes \alpha) = \\
(m_Z \otimes & 1_Z) \circ (1_Z \otimes (u_Z \circ \varepsilon_A) \otimes 1_Z)
\circ (1_Z \otimes \alpha) = 1_{Z \otimes Z}.
\end{align*}
This proves that $\kappa_r$ is bijective and 
concludes the proof of Theorem 1.2 $\square$

\bigskip

Combining Theorem 1.2 and a special case of a theorem
of P. Schauenburg \cite{[Sc2]}, we get the following result.
It would be interesting to find a direct proof.

\begin{coro}
Let $(A,B,Z,T)$ be a Hopf-Galois system. Then $A$ and $B$ are
Hopf algebras.
\end{coro}

Another theorem of P. Schauenburg (Theorem 5.5 in
\cite{[Sc1]}) ensures that if $A$ and $B$ are Hopf algebras
such that there exists an $A$-$B$-biGalois extension,
the comodule categories over $A$ and $B$ are monoidally equivalent.
This theorem, combined with Theorem 1.2, 
yields the following result.

\begin{coro}
Let $(A,B,Z,T)$ be a Hopf-Galois system. Then 
the categories $\mathrm{Comod}(A)$ and $\mathrm{Comod}(B)$ 
are monoidally equivalent.
\end{coro}

Let us now explain the reconstruction of a Hopf-Galois system
from a Galois extension. 
We use Tannaka duality techniques, for which
our references are \cite{[JS]} and \cite{[Sc0]}.
We first consider the following
general situation. Let $\mathcal C$ be a small category
and 
let $F,G : \mathcal C \longrightarrow {\rm Vect}_{\rm f}(k)$ be some functors.  
Following \cite{[JS]}, Section 3, we associate a vector space
$\mathsf{Hom}^{\vee}\!(F,G)$ to such a pair:
$$
\mathsf{Hom}^{\vee}\!(F,G) = \bigoplus_{X \in {\rm ob}(\mathcal C)} 
{\rm Hom}_k(G(X),F(X))/\mathcal N,$$
where
$\mathcal N$ be the linear subspace of $\bigoplus_{X \in {\rm ob}(\mathcal C)} {\rm
Hom}_k(G(X),F(X))$ generated by the elements
$F(f) \circ u - u \circ
G(f)$, with $f \in {\rm Hom}_{\mathcal C}(X,Y)$ and $u \in {\rm
Hom}_k(G(Y),F(X))$.
The class of an element $u \in {\rm Hom}_k(G(X), F(X))$ is denoted
by $[X,u]$ in  $\mathsf{Hom}^{\vee}\!(F,G)$.

This vector space
represents the functor ${\rm Vect}_{\rm f}(k) \longrightarrow {\rm Sets}$,
$V \longmapsto {\rm Nat}(F, G \otimes V)$
(see \cite{[JS]}). 
Now let $K : \mathcal C \longrightarrow {\rm Vect}_{\rm f}(k)$ be another
functor. The universal property of $\mathsf{Hom}^{\vee}\!(F,G)$ gives
 a linear map 
$$\delta_{F,G}^K : \mathsf{Hom}^{\vee}\!(F,G)
\longrightarrow \mathsf{Hom}^{\vee}\!(K,G) 
\otimes \mathsf{Hom}^{\vee}\!(F,K),$$
coassociative in an obvious sense. The map $\delta_{F,G}^K$
may be described as follows. Let $X$ in ob$(\mathcal C)$, let
$\phi \in G(X)^*$, let $x \in F(X)$ and let $e_1, \ldots , e_n$
be a basis of $K(X)$. Then
$$\delta_{F,G}^K([X, \phi \otimes x])=
\sum_{i=1}^n [X,\phi \otimes e_i] \otimes [X, e_i^* \otimes x].$$   
As a particular case of the previous construction,
$\mathsf{End}^{\vee}\!(F) := \mathsf{Hom}^{\vee}\!(F,F)$ is a coalgebra
(the counit is induced by the trace, see \cite{[JS]}, Section 4)

Assume now that $\mathcal C$ is a monoidal category and 
that $F$ and $G$ are monoidal functors. Then 
$\mathsf{Hom}^{\vee}\!(F,G)$ inherits an algebra structure, which
may be described by the following formula:
$$[X,u].[Y,v] = [X \otimes Y, \widetilde F_{X,Y} \circ (u \otimes v)
\circ \widetilde G_{X,Y}^{-1}],$$
where the isomorphisms $\widetilde F_{X,Y} : F(X) \otimes F(Y) 
\longrightarrow F(X \otimes Y)$ and 
$\widetilde G_{X,Y} : G(X) \otimes G(Y) 
\longrightarrow G(X \otimes Y)$ are part of the monoidal functors
$F$ and $G$. It is easy to see that the maps
$\delta_{F,G}^K$ are algebra maps, and hence
$\mathsf{End}^{\vee}\!(F)$ is a bialgebra.

Assume finally that $\mathcal C$ is a rigid monoidal category. This
means that every object $X$ has a left dual (\cite{[JS], [Ka]}), 
i.e. there exist a triplet $(X^*, e_X, d_X)$ where 
$X^* \in {\rm ob}(\mathcal C)$, while 
$e_X:  {X^* } \otimes X \longrightarrow I$ ($I$ is the monoidal
unit of $\mathcal C$) and 
$d_X : I \longrightarrow X \otimes  {X^* }$ are morphisms 
of $\mathcal C$ such that:
$$ 
(1_X \otimes e_X) \circ (d_X \otimes 1_X) = 1_X \quad 
{\rm and} \quad (e_X \otimes 1_{X^*})
 \circ (1_{X^*} \otimes d_X) = 1_{X^*}.$$
The rigidity of $\mathcal C$ allows one to define a duality
endofunctor of $\mathcal C$, which will be used in the proof of
the following result.

\begin{prop}
Let $\mathcal C$ be a rigid monoidal category and let 
$F,G : \mathcal C \longrightarrow {\rm Vect}_{\rm f}(k)$ be monoidal 
functors. Then
$(\mathsf{End}^{\vee}\!(F), \mathsf{End}^{\vee}\!(G),
\mathsf{Hom}^{\vee}\!(G,F), \mathsf{Hom}^{\vee}\!(F,G))$ is a 
Hopf-Galois system.
\end{prop}

\noindent
\textbf{Proof}. We retain the notations of Definition 1.1.
We put $\alpha := \delta_{G,F}^F$, 
$\beta :=\delta_{G,F}^G$, $\gamma := \delta_{F,F}^G$ and 
$\delta := \delta_{G,G}^F$. It is clear that the axioms
\textsf{(HG1)}-\textsf{(HG3)} are satisfied. Hence it remains to construct
the linear map $S: \mathsf{Hom}^{\vee}\!(F,G) \longrightarrow
\mathsf{Hom}^{\vee}\!(G,F)$. Let $X \in {\rm ob}(\mathcal C)$. Then
we have natural isomorphisms 
$$\lambda_X^F : F(X)^* \longrightarrow F(X^*) \quad {\rm and} \quad
\lambda_X^G : G(X)^* \longrightarrow G(X^*)$$
such that the following diagrams commute:
 $$\xymatrix{
F(X)^* \otimes F(X) \ar[r]^-{e_{F(X)}} \ar[d]^{\lambda_X^F \otimes 1_{F(X)}} &
I \ar[r]^{\widetilde F_0}  & F(I) \\
F(X^*) \otimes F(X) \ar[rr]^-{\widetilde F_{X^*,X}} && F(X^* \otimes X) \ar[u]^{F(e_X)}} \quad \quad
\xymatrix{
F(X) \otimes F(X)^*  \ar[d]^{1_{F(X)} \otimes \lambda_X^F} &
I  \ar[l]_-{d_{F(X)}} \ar[r]^{\widetilde F_0}  & F(I) \ar[d]^{F(d_X)} \\
F(X) \otimes F(X^*) \ar[rr]^-{\widetilde F_{X,X^*}} && 
F(X \otimes X^*)}
$$
Let $u \in {\rm Hom}_k(G(X), F(X))$. We put 
$$S([X,u]) = [X^*, \lambda_X^G \circ {^t\!u} \circ (\lambda_X^F)^{-1}].$$
It is easy to see that $S$ is a well defined linear map.
Now let $\phi \in F(X)^*$, let $x \in F(X)$ and let $e_1, \ldots, e_n$
be a basis of $G(X)$. Then we have
\begin{align*}
m & \circ (1 \otimes S) \circ \gamma ([X, \phi \otimes x]) 
= \sum_{i=1}^n [X, \phi \otimes e_i] [X^*, \lambda_X^G \circ
(e_i^* \otimes x) \circ (\lambda_X^F)^{-1}] = \\
= & \sum_{i=1}^n
[X \otimes X^*, \widetilde G_{X,X^*} \circ
(1_{F(X)} \otimes \lambda_X^G) \circ
((\phi \otimes e_i) \otimes (x \otimes e_i^*))
\circ (1_{F(X)} \otimes (\lambda_X^F)^{-1}) \circ
\widetilde F_{X,X^*}^{-1}] \\
= & [X \otimes X^*, \widetilde G_{X,X^*} \circ
(1_{F(X)} \otimes \lambda_X^G) \circ d_{G(X)} \circ (\phi \otimes x)
\circ (1_{F(X)} \otimes (\lambda_X^F)^{-1}) \circ
\widetilde F_{X,X^*}^{-1}] \\
= & [X \otimes X^*, G(d_X) \circ \widetilde G_0 \circ (\phi \otimes x) 
\circ (1_{F(X)} \otimes (\lambda_X^F)^{-1}) \circ
\widetilde F_{X,X^*}^{-1}] \\
= & [I, \widetilde G_0 \circ (\phi \otimes x)
\circ (1_{F(X)} \otimes (\lambda_X^F)^{-1}) \circ
\widetilde F_{X,X^*}^{-1} \circ F(d_X)] \\ 
= & [I, \widetilde G_0 \circ (\phi \otimes x) \circ d_{F(X)}
\circ \widetilde F_0^{-1}] 
= \phi(x) [I,\widetilde G_0 \circ \widetilde F_0^{-1}]
=\varepsilon([X, \phi \otimes x])1.\\
\end{align*}
Since the elements $[X, \phi \otimes x]$ linearly span 
$\mathsf{End}^{\vee}\!(F)$, we have the commutativity of the first
diagram of \textsf{HG4}. The commutativity of the second diagram is proved
similarly. $\square$

\medskip

\begin{rem}
{\rm 
Proposition 1.5 generalizes \cite{[U2]}, using exactly the same
idea. More generally, Proposition 1.5 is still valid with weaker
hypothesis on the target category (which we have assumed here to be
${\rm Vect}_{\rm f}(k)$):  see \cite{[Sc0]}. Of course the proof
is more difficult to write: see the proof of Theorem 2.4.2 in
\cite{[Sc0]}.
Our proof, using rank one operators in the case of ${\rm Vect}_{\rm f}(k)$, 
is not very elegant, but is quite straightforward.}
\end{rem}

\begin{rem}
{\rm It is easily seen that the map
$S: \mathsf{Hom}^{\vee}\!(F,G) \longrightarrow
\mathsf{Hom}^{\vee}\!(G,F)^{\rm op}$ 
constructed in the proof of Proposition 1.5 is an algebra morphism.}
\end{rem}

We can now recover a Hopf-Galois system starting from
a Galois extension. 

\begin{coro}
Let $A$ be Hopf algebra and let $Z$ be a left $A$-Galois extension.
Then there exists a Hopf algebra $B$ and an algebra $T$ such
that $(A,B,Z,T)$ is a Hopf-Galois system.
\end{coro}

\noindent
\textbf{Proof}.  
First consider the forgetful functor 
$\omega : {\rm Comod}_{\rm f}(A) \longrightarrow
{\rm Vect}_{\rm f}(k)$. By tannakian reconstruction theorems
\cite{[JS],[Sc0]} the Hopf algebras $A$ and 
$\mathsf{End}^{\vee}\!(\omega)$ are isomorphic: hence we identify
these two Hopf algebras.
Now consider the $A$-Galois extension $Z$.
Following Ulbrich \cite{[U]}, we associate a fibre functor 
$\eta_Z : {\rm Comod}_{\rm f}(A) \longrightarrow {\rm Vect}_{\rm f}(k)$
to $Z$ ($\eta_Z$ is a monoidal, $k$-linear, exact and faithful functor).
For an $A$-comodule $V$, we have $\eta_Z(V) = V \wedge Z$, where
$V \wedge Z$ is the kernel of the double arrow:
$$
\alpha_V \otimes 1_Z,  \ 1_V \otimes \alpha_Z \ : \ V \otimes Z 
\rightrightarrows
V
\otimes A \otimes Z
$$
($V\wedge Z$ is the cotensor product of \cite{[Ta2]}). We have an 
obvious monoidal natural transformation 
$\eta_Z \longrightarrow \omega \otimes Z$ and thus the universal
property of $\mathsf{Hom}^{\vee}\!(\eta_Z,\omega)$ yields an $A$-colinear
algebra morphism 
$\mathsf{Hom}^{\vee}\!(\eta_Z,\omega) \longrightarrow Z$. Since by 
Proposition 1.5 and Theorem 1.2 $\mathsf{Hom}^{\vee}\!(\eta_Z,\omega)$
is a left $A$-Galois extension, and since the category of $A$-Galois
extensions is a groupoid \cite{[U0]}, then
$\mathsf{Hom}^{\vee}\!(\eta_Z,\omega) \cong Z$.
Then, with the obvious identifications, 
$(\mathsf{End}^{\vee}\!(\omega),
\mathsf{End}^{\vee}\!(\eta_Z),
\mathsf{Hom}^{\vee}\!(\eta_Z,\omega),
\mathsf{Hom}^{\vee}\!(\omega, \eta_Z))$ is the Hopf-Galois
system we have announced. $\square$

\begin{rem}
{\rm In \cite{[Sc1]}, Schauenburg constructs the Hopf algebra $B$
(and the algebra $T$ in Section 4) using different techniques.
His techniques allow him to work with Hopf algebras over a ring
(with a faithful flatness assumption). It is certainly possible
to get the whole Hopf-Galois system using his techniques. On the other hand,
when the base ring is a field, its seems that the tannakian methods used
here are easier to use (it may be a question of personal taste).}
\end{rem}

Using Remark 1.7 and the proof of the last corollary, we easily have
the following result, generalizing the classical fact that
the antipode of a Hopf algebra is an algebra anti-morphism. 
Again it would be interesting to have a direct proof.

\begin{coro}
Let $(A,B,Z,T)$ be a Hopf-Galois system. Then $S : T \longrightarrow
Z^{\rm op}$ is an algebra morphism.
\end{coro} 

\begin{rem} 
{\rm We have done enough work to prove easily that
if $A$ and $B$ are Hopf algebras such that there exists an 
$A$-$B$-biGalois extension, then the 
comodule categories over $A$ and $B$ are monoidally equivalent.
This is the part of Schauenburg's Theorem 5.5 in \cite{[Sc1]}
that was used to prove Corollary 1.4, certainly the most important
result of the present paper. So we include a proof for the sake of 
completeness.

Let $Z$ be an $A$-$B$-biGalois extension. The fibre functor
$\eta_Z :{\rm Comod}_{\rm f}(A) \longrightarrow {\rm Vect}_{\rm f}(k)$
factorizes through ${\rm Comod}_f(B)$, and thus,
by the universal property of $\mathsf{End}^{\vee}\!(\eta_Z)$,
there exists
a Hopf algebra morphism $\phi : \mathsf{End}^{\vee}\!(\eta_Z)
\longrightarrow B$ such that$(1_Z \otimes \phi) \circ \alpha'
= \beta$, where $\alpha'$ stands for the canonical coaction of
$\mathsf{End}^{\vee}\!(\eta_Z)$ on $Z$ (recall that $\eta_Z(A) = Z$). 
Since $Z$ may be identified with $\mathsf{Hom}^{\vee}\!(\eta_Z, \omega)$
(proof of Corollary 1.8), it follows from Proposition 1.5 
and Theorem 1.2 that
$Z$ is a right $\mathsf{End}^{\vee}\!(\eta_Z)$-Galois extension.
Then we have $(1_Z \otimes \phi) \circ \kappa_r' = \kappa_r$,
($\kappa_r'$ stands for the Galois map of 
$\mathsf{End}^{\vee}\!(\eta_Z)$ relative to $Z$) and since 
$Z$ is a right $B$-Galois extension, it follows that $1_Z \otimes \phi$
is bijective, and so is $\phi : \mathsf{End}^{\vee}\!(\eta_Z) \cong B$.
Now since $\eta_Z$ is a fibre functor,
tannakian theorems \cite{[JS],[Sc0]} ensure that 
${\rm Comod}_{\rm f}(A)$ and 
${\rm Comod}_{\rm f}(\mathsf{End}^{\vee}\!(\eta_Z))$ are monoidally
equivalent. This concludes our proof
since the category ${\rm Comod}(A)$
is the category of Ind-objects of ${\rm Comod}_{\rm f}(A)$.}
\end{rem}

\section{Hopf-Galois systems and 2-Cocycles}

BiGalois extensions are associated to 2-cocycles in \cite{[Sc1]}.
We review this construction in the framework of Hopf-Galois systems.
After this, we study a concrete example for the
function algebra on the symmetric group.

Let $A$ be a Hopf algebra. We use Sweedler's notation
$\Delta(a) = a_{(1)} \otimes a_{(2)}$. Recall (see e.g. \cite{[Do]})
that a 2-cocycle is a convolution invertible linear map
$\sigma : A \otimes A \longrightarrow k$ satisfying
$$\sigma(a_{(1)}, b_{(1)}) \sigma(a_{(2)}b_{(2)},c) =
\sigma(b_{(1)},c_{(1)}) \sigma(a,b_{(2)} c_{(2)})$$
and $\sigma(a,1) = \sigma(1,a) = \varepsilon(a)$, for all $a,b,c \in A$.
The convolution inverse of $\sigma$, denoted $\bar{\sigma}$, satisfies 
$$\bar{\sigma}(a_{(1)}b_{(1)},c) \bar{\sigma}(a_{(2)},b_{(2)}) =
\bar{\sigma}(a,b_{(1)}c_{(1)}) \bar{\sigma}(b_{(2)}, c_{(2)})$$
and 
$\bar{\sigma}(a,1) = \bar{\sigma}(1,a) = \varepsilon(a)$, 
for all $a,b,c \in A$.

Following \cite{[Do]} and \cite{[Sc1]}, we associate various algebras to
a 2-cocycle. First consider the algebra 
$_{\sigma} \! A$. As a vector space $_{\sigma} \! A = A$ and the product
of $_{\sigma}A$ is defined to be
$$a _{\sigma} \! . b = \sigma(a_{(1)}, b_{(1)}) a_{(2)} b_{(2)}, 
\quad a,b \in A.$$
We also have the algebra $A_{\bar{\sigma}}$. As a vector space we have
$A_{\bar{\sigma}} = A$ and the product of 
$A_{\bar{\sigma}}$ is defined to be
$$a ._{\bar{\sigma}} b= \bar{\sigma}(a_{(2)}, b_{(2)}) a_{(1)} b_{(1)}, 
\quad a,b \in A.$$
Then $A_{\bar{\sigma}}$ is a left $A$-comodule algebra with coaction
$\alpha$ defined by $\alpha = \Delta$. Finally we have the 
Hopf algebra $_{\sigma} \! A_{\bar{\sigma}}$ (denoted
$A^{\sigma}$ in \cite{[Do]}). As a coalgebra $_{\sigma} \! A_{\bar{\sigma}}
= A$. The product of $_{\sigma} \! A_{\bar{\sigma}}$ is defined to be
$$a . b= \sigma(a_{(1)}, b_{(1)})
\bar{\sigma}(a_{(3)}, b_{(3)}) a_{(2)} b_{(2)}, 
\quad a,b \in A,$$
and we have the following formula for the antipode of 
$_{\sigma} \! A_{\bar{\sigma}}$:
$$S^{\sigma}(a) = \sigma(a_{(1)}, S(a_{(2)}))
\bar{\sigma}(S(a_{(4)}), a_{(5)}) S(a_{(3)}).$$
The algebra $A_{\bar{\sigma}}$ is a right 
$_{\sigma} \! A_{\bar{\sigma}}$-comodule algebra, with coaction
defined by $\beta = \Delta$. In this way
$A_{\bar{\sigma}}$ is an $A$-$_{\sigma} \! A_{\bar{\sigma}}$-bicomodule 
algebra. We have the following result.

\begin{prop}
Let $A$ be Hopf algebra and let $\sigma : A \otimes A \longrightarrow k$
be a 2-cocycle. Then $(A, \ _{\sigma} \! A_{\bar{\sigma}}, A_{\bar{\sigma}}, \
_{\sigma} \! A)$ is a Hopf-Galois system.
\end{prop}

\noindent
\textbf{Proof}. We put $\gamma = \delta = \Delta$. It is easy to see
that the axiom \textsf{HG3} is satisfied. Now define a linear map
$\phi : \  _{\sigma} \! A \longrightarrow A_{\bar{\sigma}}$
by $\phi(a) = \sigma(a_{(1)}, S(a_{(2)})) S(a_{(3)})$, for $a\in A$.
Then 
\begin{align*}
m_{A_{\bar{\sigma}}} &  
\circ (1_{A_{\bar{\sigma}}} \otimes \phi) \circ \gamma(a)
= a_{(1)} ._{\bar{\sigma}} \phi(a_{(2)})
= a_{(1)} ._{\bar{\sigma}} \sigma(a_{(2)}, S(a_{(3)})) S(a_{(4)}) = \\
\sigma & (a_{(3)}, S(a_{(4)}))
\bar{\sigma}(a_{(2)},S(a_{(5)})) a_{(1)} S(a_{(6)}) = \\
\bar{\sigma} * & \sigma (a_{(2)},S(a_{(3)}))a_{(1)} S(a_{(4)}) =
a_{(1)}S(a_{(2)}) = \varepsilon(a)1.
\end{align*}
We also have
\begin{align*}
m_{A_{\bar{\sigma}}} &  
\circ (\phi \otimes 1_{A_{\bar{\sigma}}}) \circ \delta(a)
= \phi(a_{(1)})._{\bar{\sigma}} a_{(2)} = 
\sigma  (a_{(1)}, S(a_{(2)})) S(a_{(3)}) ._{\bar{\sigma}} a_{(4)} = \\
\sigma  & (a_{(1)}, S(a_{(2)})) 
\bar{\sigma}(S(a_{(3)}), a_{(6)}) S(a_{(4)}) a_{(5)} =
 \sigma(a_{(1)}, S(a_{(2)})) 
\bar{\sigma}(S(a_{(3)}), a_{(4)}) = \varepsilon(a)1,
\end{align*}
by (a5) of Theorem 1.6 in \cite{[Do]}. Thus  
$(A, \ _{\sigma} \! A_{\bar{\sigma}}, A_{\bar{\sigma}}, \
_{\sigma} \! A)$ is a Hopf-Galois system. $\square$

\medskip

Let now study an explicit example. In fact the cocycle 
will only be used when proving that a certain algebra is non-zero.
In \cite{[Bi1]} we have constructed 2-cocycle deformations of the
function algebra on the symmetric group. We generalize these
results here.

Let us fix some notations. Until the end of the section $k$
will be a characteristic zero field.
We fix $m,n \in \mathbb Z^*$ with $m,n \geq 2$ and a primitive
$m$-th root of unity $\xi$ contained in $k$. 
We will work with the symmetric group
$S_{mn}$. For a real number $x$, we put $E^+(x) = n$ where
$n \in \mathbb Z$ is such that $x \in ]n-1,n]$. For 
$i \in \{1, \ldots , mn \}$, we put $i^* := E^{+}(\frac{i}{m})
\in \{1, \ldots , n \}$. 

We say that a matrix ${\bf p}=(p_{ij}) \in M_n(k)$ is an AST-matrix
(after Artin-Schelter-Tate \cite{[AST]}) if 
$p_{ii} = 1$ and $p_{ij}p_{ji} = 1$ for all $i$ and $j$. An AST-matrix
is said to be of order $m$ if $p_{ij}^m = 1$ for all $i$ and $j$.
The trivial AST-matrix (i.e. $p_{ij} =1$ for all $i$ and $j$)
is denoted by ${\bf 1}$.

Let ${\bf p} \in M_n(k)$ be an AST-matrix of order $m$. Let $i,j,k,l \in
\{1, \ldots , mn \}$. We put
$$R_{ij}^{lk}({\bf p}) := \delta_{i^*k^*} \delta_{j^*l^*}
\sum_{r,s = 0}^{m-1} \xi^{r(i-k) + s(j-l)} p_{j^*i^*}^{rs}.$$

\begin{defi}
Let ${\bf p},{\bf q} \in M_n(k)$ be AST matrices of order $m$. The algebra
$\mathcal O_{{\bf q}, {\bf p}}(S_{mn})$ is defined to be the universal
algebra with generators $(x_{ij})_{1 \leq i,j \leq mn}$ and satisfying the
relations:  
\begin{align}
x_{ij}   x_{ik} = \delta_{jk} & x_{ij} \quad ; \quad
x_{ji} x_{ki} = \delta_{jk} x_{ji} \quad ; \quad
\sum_{l=1}^{mn} x_{il} = 1 = \sum_{l=1}^{mn} x_{li} \quad , \
1 \leq i,j,k \leq n. \\
\sum_{k,l} & R_{ij}^{lk}({\bf p}) x_{\alpha l} x_{\beta k} =
\sum_{k,l} R_{lk}^{\alpha \beta}({\bf q}) x_{li} x_{kj} \ , \
1 \leq i,j, \alpha , \beta \leq n.
\end{align}
\end{defi}

When  ${\bf p} = {\bf q}$, then it is easily seen that 
$\mathcal O_{{\bf p}}(S_{mn})
:=\mathcal O_{{\bf p}, {\bf p}}(S_{mn})$ is a Hopf algebra,
with coproduct defined by
$\Delta(x_{ij}) = \sum_k x_{ik} \otimes x_{kj}$, counit
defined by 
$\varepsilon(x_{ij}) = \delta_{ij}$ and antipode 
defined by $S(x_{ij})= x_{ji}$
(note that $R_{ij}^{lk}({\bf p}) = R_{kl}^{ji}({\bf p})$).
Note that the relations (2) are just FRT relations \cite{[RTF]}. 
If $m = 2$, the present Hopf algebras coincide with the Hopf
algebras $\mathcal O_{{\bf p}}(S_{2n})$ of \cite{[Bi1]}.
The algebras $\mathcal O_{{\bf q}, {\bf p}}(S_{mn})$ will be shown
to part of a Hopf-Galois system. Before we need
a lemma.

\begin{lemm}
Let ${\bf p} \in M_n(k)$ be an AST matrix of order $m$.
Then $\mathcal O_{{\bf p}, {\bf 1}}(S_{mn})$ is a non-zero
algebra.
\end{lemm}

\noindent
\textbf{Proof}. We will use an appropriate 2-cocycle.
For $1 \leq i \leq n$, put $t_i = (m(i-1) +1, \ldots , mi) \in S_{mn}$, and
let $H = \langle t_1 , \ldots , t_n \rangle$ ($H \cong 
(\mathbb Z / m \mathbb Z)^n$). 
Following  Artin-Schelter-Tate \cite{[AST]}, we define
$\sigma_{\bf p} : H \times H \longrightarrow k^*$ to be the
unique bimultiplicative map such that
$\sigma_{\bf p}(t_i,t_j) = p_{ij}$ for $i<j$ and
$\sigma_{\bf p}(t_i,t_j) = 1$ for $i \geq j$.
Now consider the surjective Hopf algebra morphism
$$\pi :\mathcal O(S_{mn}) \longrightarrow k \lbrack H \rbrack,
\quad x_{ij} \longmapsto 
\frac{\delta_{i^*j^*}}{m}\sum_{k=0}^{m-1}
\xi^{k(j-i)}t_{i*}^k \ , \ 1 \leq i,j,k,l \leq mn,$$
where $\mathcal O(S_{mn}) = \mathcal O_{\bf 1}(S_{mn})$ is the 
function algebra on the symmetric group.
Composing now $\pi \otimes \pi$ with the unique
$k$-linear extension of $\sigma_{\bf p}$ to
$k[H] \otimes k[H]$, we get a 2-cocycle on 
$\mathcal O(S_{mn})$, still denoted $\sigma_{\bf p}$. 
This is the method of construction of 2-cocycles induced
by abelian subgroups of Enock-Vainerman \cite{[EV]}.
We have
$$
\sigma_{\bf p}(x_{ij}, x_{kl}) = \delta_{ij} \delta_{kl} \ {\rm if} \
i^* \geq k^* \ {\rm and} \
\sigma_{\bf p}(x_{ij}, x_{kl}) =
\frac{\delta_{i^*j^*} \delta_{k^*l^*}}{m^2}\sum_{r,s=0}^{m-1}
\xi^{r(j-i) + s(l-k)}p_{i^*k^*}^{rs} \ {\rm if} \ i^* < k^*.$$
It is then a straightforward but tedious computation
to check that the generators of
$_{\sigma_{\bf p}}\mathcal O(S_{mn})$ satisfy the defining relations of
$\mathcal O_{{\bf p}, {\bf 1}}(S_{mn})$, and thus 
$\mathcal O_{{\bf p}, {\bf 1}}(S_{mn})$
is a non-zero algebra. $\square$

\begin{prop}
Consider 
${\bf p},{\bf q} \in M_n(k)$ some AST matrices of order $m$.
Then \break $(\mathcal O_{{\bf q}}(S_{mn}),
\mathcal O_{{\bf p}}(S_{mn}),
\mathcal O_{{\bf q}, {\bf p}}(S_{mn}),
\mathcal O_{{\bf p}, {\bf q}}(S_{mn}))$ is a Hopf-Galois system.
\end{prop}

\noindent
\textbf{Proof}.
Let ${\bf r} \in M_n(k)$ be another AST matrix. It is straightforward
to check that there exist a unique algebra morphism
$$\delta_{{\bf q},{\bf p}}^{\bf r}
: \mathcal O_{{\bf q}, {\bf p}}(S_{mn})
\longrightarrow \mathcal O_{{\bf q}, {\bf r}}(S_{mn}) \otimes
\mathcal O_{{\bf r}, {\bf p}}(S_{mn})$$
such that $\delta_{{\bf q},{\bf p}}^{\bf r}(x_{ij})
= \sum_k x_{ik} \otimes x_{kj}$. 
Similarly it easy to see (using $R_{ij}^{lk}({\bf p}) = R_{kl}^{ji}({\bf p})$) 
that there exist a unique algebra isomorphism
$$\phi : \mathcal O_{{\bf p}, {\bf q}}(S_{mn}) \longrightarrow
\mathcal O_{{\bf q}, {\bf p}}(S_{mn})^{\rm op}$$
such that $\phi(x_{ij}) = x_{ji}$.
Now using $\delta_{{\bf q},{\bf p}}^{\bf 1}$, Lemma 2.3 and $\phi$, 
we see that $\mathcal O_{{\bf q}, {\bf p}}(S_{mn})$ 
and $\mathcal O_{{\bf p}, {\bf q}}(S_{mn})$
are  non-zero algebras.
We have defined all the necessary structural morphisms, and it is 
immediate to check that the axioms of a Hopf-Galois system
are satisfied. $\square$

\bigskip

Combining Proposition 2.4 and Corollary 1.4, we get the following
result.

\begin{coro}
Let ${\bf p} \in M_n(k)$ be an AST matrix of order $m$. 
Then the the category
of $\mathcal O_{\bf p}(S_{mn})$-comodules is monoidally
equivalent to the
representation category of the \break symmetric group $S_{mn}$.
\end{coro}

\section{Hopf-Galois systems for Hopf algebras of bilinear forms}

In \cite{[Bi3]} we constructed Hopf biGalois extensions for the
universal Hopf algebras associated to  non-degenerate bilinear
forms. We reconsider this construction at the Hopf-Galois system level:
this makes the considerations of \cite{[Bi3]} more transparent.
Note that in general, the Hopf-Galois systems we have here
cannot be obtained using 2-cocycles.

\medskip

Let $E \in GL_m(k)$ and let $F \in GL_n(k)$. Recall \cite{[Bi3]} that
the algebra $\mathcal B(E,F)$ is the universal algebra with generators
$x_{ij}$, $1\leq i \leq m, 1\leq j \leq n$,
and satisfying the relations
$$F^{-\!1} {^t \! x} E x = I_n \ ; \   x F^{-\!1} {^t \! x} E = I_m,$$
where $x$ is the matrix $(x_{ij})$ and $I_m$ and $I_n$ are the identity
matrices of size $m$ and $n$ respectively. For $E = F$ we have the
Hopf algebra $\mathcal B(E)$ of M. Dubois-Violette and G. Launer \cite{[DVL]}.

\begin{prop}
Let $E \in GL_m(k)$ and let $F \in GL_n(k)$ ($m,n \geq 2$) be
such that  
${\rm tr}(E ^t \!  E^{-1}) \break = {\rm tr}(F ^t \!  F^{-1})$. Then
$(\mathcal B(E), \mathcal B(F), \mathcal B(E,F), \mathcal B(F,E))$
is a Hopf-Galois system.
\end{prop}

\noindent
\textbf{Proof}.
First the end of Section 4 in \cite{[Bi3]} ensures that
$\B (E,F)$ is a non-zero algebra.
Let $G \in GL_p(k)$. It is a direct computation to check
that there exists a unique algebra morphism
$\delta_{E,F}^G : \B(E,F) \longrightarrow
\B(E,G) \otimes \B(G,F)$ such that 
$\delta_{E,F}^G (x_{ij}) = \sum_{k=1}^p x_{ik} \otimes x_{kj}$, 
$1\leq i \leq m, 1\leq j \leq n$. Also there exists a unique algebra
isomorphism $\phi : \B(F,E) \longrightarrow \B(E,F)^{\rm op}$ such that
$\phi(x) = F^{-1} {^t \! x}E$.
In this way we have all the necessary structural maps
and it is immediate to check that we indeed have a Hopf-Galois system.
$\square$  

\bigskip

Using Proposition 3.1 and Corollary 1.4, we have the following
result from \cite{[Bi3]}:

\begin{coro}
1) Let $E \in GL_m(k)$ and let $F \in GL_n(k)$ ($m,n \geq 2$) be such that  
${\rm tr}(E ^t \!  E^{-1}) \break = {\rm tr}(F ^t \!  F^{-1})$. Then
the categories 
${\rm Comod}(\B(E))$ and ${\rm Comod}(\B(F))$ are monoidally 
\break equivalent.

\noindent
2) Assume that $k$ is algebraically closed.
Let $E \in GL_m(k)$ ($m \geq 2$) and let $q \in k^*$ be such
that $q^2 + {\rm tr}(E ^t \!  E^{-1})q +1 = 0$. Then
the categories 
${\rm Comod}(\B(E))$ and ${\rm Comod}(\qsl)$ are monoidally equivalent.
\end{coro}

\section{Hopf-Galois systems for cosovereign Hopf algebras}

Recall \cite{[Bi2]} that a Hopf algebra $A$ is said to be cosovereign
if there exists a character $\Phi \in A^*$ such that
$S^2 = \Phi * {\rm id} * \Phi^{-1}$. The universal (or free)
cosovereign Hopf algebras were constructed in \cite{[Bi2]}. 
We describe some of the Hopf-Galois systems associated
with this class of Hopf algebras. Our constructions are 
certainly incomplete, but they nevertheless enable us to improve
on certain known results \cite{[Ba]} on the corepresentation theory
of the universal cosovereign Hopf algebras (when $k = \C$).

\begin{defi}
Let $E \in GL_m(k)$ and let $F \in GL_n(k)$.
The algebra $H(E,F)$ is defined to be the universal algebra with generators
$u_{ij}$, $ v_{ij}$, $1\leq i \leq m, 1\leq j \leq n$,
and satisfying the relations
$$u {^t \! v} = I_m = v F {^t \! u} E^{-1} \quad ; \quad
{^t \! v}u = I_n = F {^t \! u} E^{-1} v.$$
\end{defi}

When $E=F$, we just have the universal cosovereign 
Hopf algebras $H(F)$ of \cite{[Bi2]}.
It is known (see Proposition 3.3 in \cite{[Bi2]}) that the 
Hopf algebra $H(F)$ remains unchanged, up to isomorphism,
if the matrix $F$ is multiplied by a non-zero scalar or is replaced
by a conjugate matrix. Similarly, we have the following result.

\begin{prop}
Let $\lambda \in k^*$, let $E,P \in GL_m(k)$ and let $F,Q \in GL_n(k)$. 
Then $H(\lambda E, \lambda F) = H(E,F)$, and we have algebra isomorphisms  
$H(E,F) \cong H(PEP^{-1},QFQ^{-1})$ 
and \break $H(E,F) \cong H({^t \! E}^{-1}, {^t \! F}^{-1})$.
\end{prop}

\noindent
\textbf{Proof}.
The first assertion is obvious. It is easily seen that
there exists a unique algebra isomorphism
$f : H(E,F) \longrightarrow H(PEP^{-1},QFQ^{-1})$ such that
$f(u) = {^t \! P} u {^t \! Q}^{-1}$ and $f(v) = P^{-1} v Q$.
Also we have an algebra isomorphism
$g :  H(E,F) \longrightarrow H({^t \! E}^{-1}, {^t \! F}^{-1})$ such that
$g(u) = v$ and $g(v) = E u F^{-1}$. $\square$

\bigskip

The Hopf algebra structure of $H(F)$ is a particular
case of the following result.

\begin{prop} 
Let $E \in GL_m(k)$ and let $F \in GL_n(k)$.
Assume that $H(E,F) \not = 0$.
Then $(H(E),H(F),H(E,F),H(F,E))$ is a Hopf-Galois system.
\end{prop}

\noindent
\textbf{Proof}.
Let $G \in GL_p(k)$. Then it is easy to check that
there exists a unique algebra morphism
$\delta_{E,F}^G : H(E,F) \longrightarrow H(E,G) \otimes H(G,F)$ such
that $\delta_{E,F}^G(u_{ij}) =  \sum_{k=1}^p
u_{ik} \otimes u_{kj}$ and
$\delta_{E,F}^G(v_{ij}) =  \sum_{k=1}^p
v_{ik} \otimes v_{kj}$,
$1\leq i \leq m$, $1\leq j \leq n$. Also
there is a unique algebra morphism 
$\phi : H(F,E) \longrightarrow H(E,F)^{\rm op}$ such that
$\phi(u) = {^t \! v} $ and $\phi(v) = F {^t \! u} E^{-1}$.
Thus with the obvious structural morphisms, we have a Hopf-Galois system.
$\square$

\bigskip

Of course this last result is useful only when one knows that
$H(E,F)$ is non-zero. In view of the results of the preceding section,
it is quite natural to think that $H(E,F)$ will be a non-zero algebra
when ${\rm tr}(E) = {\rm tr}(F)$ and 
${\rm tr}(E^{-1}) = {\rm tr}(F^{-1})$. This problem will be studied
elsewhere. There is already an interesting case where we can
prove that $H(E,F) \not = 0$.

\begin{prop}
Let $E \in GL_m(k)$ and let $F \in GL_n(k)$
be such that ${\rm tr}(E) = {\rm tr}(F)$. Assume that
there exists $G \in GL_m(k)$ and $K \in GL_n(k)$
such that $E = {^t \! G} G^{-1}$
and $F = {^t \! K} K^{-1}$. Then
$H(E,F)$ is a non-zero algebra.
\end{prop}

\noindent
\textbf{Proof}.
It is easy to check that there exists a unique algebra morphism
$f : H(E,F) \longrightarrow \B(G,K)$ such that
$f(u) = x$ and $f(v) = {^t \! G} x {^t \! K}^{-1}$. 
We have ${\rm tr}(E) = {\rm tr}({^t \! G} G^{-1}) = {\rm tr}(F)
= {\rm tr}({^t \! K} K^{-1})$, so 
by \cite{[Bi3]}, we know  that
$\B(G,K)$ is a non-zero algebra. Since $f$ is surjective,
it is clear that $H(E,F)$ is a non-zero algebra. $\square$

\bigskip

Let $q \in k^*$. In the next result, 
we consider the matrix 
$F_q = \left(\begin{array}{cc} q & 0 \\
                          0 & q^{-1}\\
       \end{array} \right) \in GL_2(k)$,
and we put $H_q := H(F_q)$.

\begin{coro}
1) Let $E \in GL_m(k)$ and let $F \in GL_n(k)$
be such that ${\rm tr}(E) = {\rm tr}(F)$. Assume that
there exists $G \in GL_m(k)$ and $K \in GL_n(k)$
such that $E = {^t \! G} G^{-1}$
and $F = {^t \! K} K^{-1}$. Then the categories
${\rm Comod}(H(E))$ and ${\rm Comod}(H(F))$ are monoidally 
equivalent.

\noindent
2) Let $F \in GL_n(k)$.
Assume that $k$ is algebraically closed and that
there exists $K \in GL_n(k)$  such that $F = {^t \! K} K^{-1}$. 
Let $q \in k^*$ be such that $q^2 - {\rm tr}(F)q + 1 = 0$. Then 
the categories
${\rm Comod}(H(F))$ and ${\rm Comod}(H_q)$ are monoidally 
equivalent.
\end{coro}

\noindent
\textbf{Proof}. The first assertion follows from
Propositions 4.3-4.4 and Corollary 1.4. We have $F_q = {^t \! G} G^{-1}$
for the matrix $G = \left(\begin{array}{cc} 0 & 1 \\
                          q & 0\\
       \end{array} \right)$,
and hence the second assertion follows from the fist one. $\square$

\bigskip

Until the end of the section, we assume that $k =  \C$.
Let us recall that a Hopf $*$-algebra is a Hopf algebra $A$, which
is also a $*$-algebra and such that the comultiplication is a 
$*$-homomorphism. 
Recall \cite{[KS]} that a Hopf $*$-algebra $A$ is said to be a CQG
algebra if for every finite-dimensional $A$-comodule with associate
matrix of coefficients $a \in M_n(A)$, there exists $K \in GL_n(\C)$
such that the matrix $KaK^{-1}$ is unitary.
A CQG algebra may be seen as the algebra of representative functions 
on a compact quantum group.  

Let $F \in GL_n(\C)$. We have seen in \cite{[Bi2]} (Proposition 3.6)
that $H(F)$ admits
a CQG algebra structure if and only if $F$ is conjugate to
a relatively positive matrix (a matrix $M$ is said to be relatively positive
if there exists $\lambda \in \C^*$ such that $\lambda M$ is a positive matrix).
In this case $H(F)$ is the dense Hopf $*$-algebra of one the
universal compact quantum groups introduced by A. Van Daele and 
S. Wang \cite{[VDW]}, and the corepresentation theory has been
worked out by T. Banica \cite{[Ba]}: the irreducible comodules
are labelled by the free product $\N * \N$.
We can combine Banica's results \cite{[Ba]}
and  Corollary 4.5 to get the cosemisimplicity
of some universal cosovereign Hopf algebras which do not admit
a CQG algebra structure, as well as their corepresentation theory.

\begin{ex}
{\rm 
Let $q , \alpha \in \mathbb C ^*$
Consider the matrix 
$F = \left(\begin{array}{cc} q & \alpha \\
                          0 & q^{-1} \\
       \end{array} \right)$.
Since $H(F) = H(-F)$, we can assume that $q \not = 1$ without changing the 
Hopf algebra $H(F)$. The matrix $F$ is not relatively positive, but 
satisfies the condition of Corollary 4.5, for 
$K = \left(\begin{array}{cc} \alpha \frac{q}{1-q} & 1 \\
                          q & 0\\
       \end{array} \right)$,
and hence ${\rm Comod}(H(F)) \cong^{\otimes} {\rm Comod}(H_q)$.
If $q \in \mathbb R^*$, then
$H_q$ is a CQG algebra and we can use the results of \cite{[Ba]}. 

Another example is constructed as follows. 
Let $\xi$ be a primitive $m$-th root of unity, \break
$m \geq 5$. Consider the diagonal matrix 
$F = {\rm Diag} (\xi, 1, \xi^{-1})$. 
Then $F$ is not a relatively positive matrix,
but $F$ satisfies the condition 
of Corollary 4.5 (easy to check) and the solutions of
$q^2 - (1+ \xi + \xi^{-1})q +1 = 0$ are real numbers. Hence 
we have ${\rm Comod}(H(F)) \cong^{\otimes} {\rm Comod}(H_q)$, and  
$H_q$ is CQG algera since $q$ is a real number: we can use the results
of \cite{[Ba]}.}
\end{ex}

\section{Hopf-Galois systems for free Hopf algebras}

M. Takeuchi has constructed in \cite{[T]} the free Hopf algebra
generated by a coalgebra. We consider here the case of a 
matrix coalgebra $M_m(k)^*$, and construct the corresponding Hopf-Galois
system.

\begin{defi}
Let $m,n \in \N^*$. The algebra $H(m,n)$ is defined to be the universal 
algebra with generators $x_{ij}^{(\alpha)}$, 
$1\leq i \leq m, 1\leq j \leq n$, $\alpha \in \N$, 
and submitted to the relations:
$$x^{(\alpha)} {^t \! x^{(\alpha +1)}} = I_m \ , \ 
{^t \! x^{(\alpha +1)}} x^{(\alpha)} = I_n \ , \ \alpha \in \N.$$
\end{defi}

When $m=n$, we have the free Hopf algebra $H(m,m) = 
H(m) = H(M_m(k)^*)$ of \cite{[T]}.
This Hopf algebra is also considered in \cite{[VDW]}, under a different
notation. See \cite{[T]} or \cite{[VDW]} for the structural morphisms 
of the Hopf algebra $H(m)$. In fact we have the following
more general result.

\begin{prop}
Let $m,n \geq 2$. Then 
$(H(m), H(n), H(m,n), H(n,m))$ is a Hopf-Galois system.
\end{prop}

\noindent
\textbf{Proof}. Let us first check $H(m,n)$ is a non-zero algebra.
Let $E \in GL_m(k)$ and $F \in GL_n(k)$ be such that
${\rm tr}(E ^t \!  E^{-1}) = {\rm tr}(F ^t \!  F^{-1})$. It is a direct
computation to check that there exists a unique algebra morphism
$f : H(m,n) \longrightarrow \B(E,F)$ such that
$$f(x^{(2k)}) = (E^{-1} {^t \! E})^k x (F^{-1} {^t \! F})^k
\ {\rm and} \
f(x^{(2k+1)}) =  {^t \! E} (E^{-1} {^t \! E})^k x (F^{-1} {^t \! F})^k
{^t \! F^{-1}}, \ k \in \N.$$
Thus, since $f$ is surjective and $\B(E,F)$ is a non-zero algebra 
\cite{[Bi3]},
it is clear that $H(m,n)$ is a non-zero algebra.
Let $p \geq 2$. There is a unique algebra morphism
$\delta_{m,n}^p : H(m,n) \longrightarrow H(m,p) \otimes H(p,n)$ such
that $\delta_{m,n}^p(x_{ij}^{(\alpha)}) =  \sum_{k=1}^p
x_{ik}^{(\alpha)} \otimes x_{kj}^{(\alpha)}$,
$1\leq i \leq m, 1\leq j \leq n$, $\alpha \in \N$. Also
there is a unique algebra morphism 
$\phi : H(n,m) \longrightarrow H(m,n)^{\rm op}$ such that
$\phi(x^{(\alpha)}) = {^t \! x^{(\alpha +1)}}$.
Thus with the obvious structural morphisms, we have a Hopf-Galois system.
$\square$

\bigskip

Combining Proposition 5.2 and Corollary 1.4, we have:

\begin{coro}
Let $m \geq 2$. Then the categories   
${\rm Comod}(H(m))$ and ${\rm Comod}(H(2))$ are monoidally 
equivalent.
\end{coro}

There is also a version of free Hopf algebras with a bijective antipode,
considered in \cite{[Ma]} and \cite{[VDW]}. It is left as an
exercice to the reader, using the preceding techniques,
 to construct the corresponding Hopf-Galois systems.

\end{document}